\documentclass[12pt]{amsart}
\usepackage[english]{babel}
\usepackage{graphicx}
\usepackage{amscd}
\usepackage{amsfonts}
\usepackage{latexsym}
\usepackage{amsthm, amssymb,amsmath}
\usepackage[all]{xy}
\usepackage{enumerate}
\usepackage{tikz}
\usepackage{hyperref}

\newcommand{\Q}{{\mathbb Q}}
\newcommand{\Z}{{\mathbb Z}}
\newcommand{\F}{{\mathbb F}}
\newcommand{\R}{{\mathbb R}}
\newcommand{\C}{{\mathbb C}}
\def\O{{\ensuremath{\mathcal{O} }}}
\DeclareMathOperator{\Gal}{Gal}
\DeclareMathOperator{\Ind}{Ind}
\DeclareMathOperator{\cnd}{N}
\DeclareMathOperator{\GL}{GL}
\DeclareMathOperator{\coh}{H}
\DeclareMathOperator{\Frob}{Frob}
\newcommand{\abcd}[4]{\ensuremath{\begin{pmatrix}  #1&#2\\#3&#4\end{pmatrix}}}
\newcommand{\abc}[4]{\ensuremath{\left(\begin{smallmatrix}  #1&#2\\#3&#4\end{smallmatrix}\right)}}

\def\p{{\ensuremath{\mathfrak{p}}}}

\DeclareMathOperator{\nm}{Nm}
\DeclareMathOperator{\Id}{Id}
\DeclareMathOperator{\Tr}{Trace}
\def\triv{{\mathbf{1}}}

\newtheorem{thm}{Theorem}

\newtheorem{cor}[thm]{Corollary}
\newtheorem*{conj}{Conjecture}
\theoremstyle{definition}
\newtheorem{ex}{Exercise}
\newtheorem{dfn}[thm]{Definition}
\newtheorem{exa}[thm]{Example}
\newtheorem*{rmk}{Remark}
\newtheorem*{rmks}{Remarks}

\makeatletter
\def\@tocline#1#2#3#4#5#6#7{\relax
  \ifnum #1>\c@tocdepth 
  \else
    \par \addpenalty\@secpenalty\addvspace{#2}%
    \begingroup \hyphenpenalty\@M
    \@ifempty{#4}{%
      \@tempdima\csname r@tocindent\number#1\endcsname\relax
    }{%
      \@tempdima#4\relax
    }%
    \parindent\z@ \leftskip#3\relax \advance\leftskip\@tempdima\relax
    \rightskip\@pnumwidth plus4em \parfillskip-\@pnumwidth
    #5\leavevmode\hskip-\@tempdima
      \ifcase #1
       \or\or \hskip 1em \or \hskip 2em \else \hskip 3em \fi%
      #6\nobreak\relax
    \dotfill\hbox to\@pnumwidth{\@tocpagenum{#7}}\par
    \nobreak
    \endgroup
  \fi}
\makeatother

\begin{document}
\begin{titlepage}
\renewcommand{\thefootnote}{\fnsymbol{footnote}}
\begin{center}
\textmd{Advanced School and Workshop on $L$\--functions and modular forms\\ICTP Trieste, September 2014}
\vfill
{\huge \textbf{$l$-adic representations and their associated invariants} }\\[0.4cm]
\begin{center}
\vfill
By\\
{\large{Vladimir Dokchitser}}\parindent=0pt\small{\footnote{University of Warwick, supported by a Royal Society University Research Fellowship,\\\url{ V.Dokchitser@warwick.ac.uk} }}
\end{center}

\begin{center}
\textmd{Notes by Samuele Anni\parindent=0pt\footnote{University of Warwick, supported by EPSRC Programme Grant ``LMF: $L$\--Functions and Modular Forms", \url{ S.Anni@warwick.ac.uk}}}\\
\vfill
\end{center}

\end{center}

\end{titlepage}

\newpage

\subsection*{Introduction }$\,$\\

These notes are from a 3-lecture course given by the author at the ICTP in Trieste, Italy, 1st--5th of September 2014, as part of a graduate summer school on ``$L$-functions and modular forms''.

The course is meant to serve as an introduction to $l$-adic Galois representations over local fields with ``$l\neq p$'', and, as the school was followed by a computational workshop, has a somewhat computational bent. It is worth mentioning that the course is not about varieties and their \'etale cohomology, but merely about the representation theory. It should also be said that nothing here is meant to be new --- in fact, it's all pretty old by now, except for perhaps a couple tricks illustrated in the exercises, that were inspired by more recent developments. 


The prerequisites for the course are fairly modest: the theory of local fields and representation theory of finite groups, such as might be covered in a Masters level course, and, for the sake of examples and motivation, the theory of elliptic curves over local fields (Silverman's book is more than sufficient). Some previous exposure to $L$-functions is desirable, as, after all, this was the main topic of the school. The ``Preliminary exercises'' in the beginning are meant to set the right level, as well as serve as a warm-up for the course.

\bigskip\bigskip\bigskip
\setcounter{tocdepth}{1}
\tableofcontents

\newpage

\section*{Preliminary exercises}

\bigskip

\subsection*{Representation theory}
\begin{ex}
Write down the irreducible representations of $S_3$ and $D_{10}$ (throughout $S_n$ (resp. $A_n$) will denote the symmetric (resp. alternating) group on $n$ elements, and $D_{2n}$ will denote the dihedral group of order $2n$).
\end{ex}
\begin{ex}
Using the fact that $S_4/V_4\cong S_3$, find the irreducible representations of $S_4$ (or their characters), where $V_4$ denotes the Klein group.
\end{ex}
\begin{ex}
Compute the induced characters of the two irreducible $2$\--dimensional representations of $D_{10}$ to $A_5$. Compute the character table of $A_5$.
\end{ex}

\bigskip\bigskip\bigskip

\subsection*{Local fields}
\setcounter{ex}{0}
\begin{ex}
Determine the degrees, ramification degrees and residue degrees of 
$$\Q_2(\zeta_3, \sqrt[3]{3})/\Q_2, \quad \quad \Q_3(\zeta_3, \sqrt[3]{3})/\Q_3\quad\mbox{ and }\quad \Q_7(\zeta_3, \sqrt[3]{3})/\Q_7.$$
Here $\zeta_3$ denotes a primitive cube root of 1.
\end{ex}
\begin{ex}
Compute the corresponding Galois groups, inertia groups, wild inertia groups and Frobenius elements.
\end{ex}

\bigskip\bigskip\bigskip

\subsection*{Elliptic curves}
\setcounter{ex}{0}
\begin{ex}
Show that the elliptic curve  $E: y^2+y=x^3-x^2$ has split multiplicative reduction at $11$ and good reduction at all other primes. Compute $\#\tilde{E}(\F_2)$, $\#\tilde{E}(\F_5)$ and $\#\tilde{E}(\F_{11})$.
\end{ex}
\begin{ex}
Find an elliptic curve with additive potentially good reduction at $5$, and another with additive potentially multiplicative reduction at $5$.
\end{ex}

\newpage 
\section{Artin representations}

Throughout these lectures $K$ will denote a field of characteristic zero.

We will abuse notation for representations and use the same symbol both for the homomorphism $G\to \GL_d(\C)$ and for the underlying vector space.

\begin{minipage}{0.7\linewidth} 
\begin{dfn}An \emph{Artin representation} $\rho$ over $K$ is a finite dimensional complex representation of $\Gal(\overline{K}/K)$ that factors through a finite quotient by an open subgroup, i.e.\ there is a finite Galois extension $F/K$ such that $\Gal(\overline{K}/F)\subseteq \ker\rho$, so that $\rho$ corresponds to a representation of  $\Gal(F/K)$.
\end{dfn}

\end{minipage}
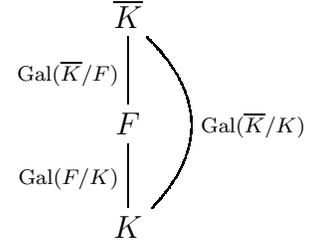
\begin{minipage}{0.3\linewidth}
$$\xymatrix{
\overline{K}\ar^{\Gal(\overline{K}/K)}@/^2pc/@{-}[dd]\ar_{\Gal(\overline{K}/F)}@{-}[d]\\
F\ar_{\Gal(F/K)}@{-}[d]\\
K}$$
\end{minipage}

\begin{exa}
$K=\Q_5$, $F=\Q_5(\zeta_3, \sqrt[3]{5})$, $$\Gal(F/K)\cong S_3 =\langle s,t|s^3=t^2=id, tst=s^{-1}\rangle$$
\begin{itemize}
\item Let $\rho$ be the two dimensional irreducible representation of $S_3$:

\begin{minipage}{0.8\linewidth} 
$$\rho(s)=\abcd{-{1}/{2}}{{{-}\sqrt{3}}/{2}}{{\sqrt{3}}/{2}}{-{1}/{2}}\quad \rho(t)=\abcd{1}{0}{0}{{-}1}.$$
\end{minipage}
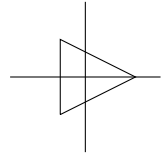
\begin{minipage}{0.2\linewidth}
\begin{tikzpicture}
\draw (-1,0) -- (1,0);
\draw (0,-1) -- (0,1);
\draw (-1/3,-1/2) -- (-1/3,1/2)--(2/3,0)--(-1/3,-1/2);
\end{tikzpicture}
\end{minipage}
It gives an irreducible representation of $\Gal(\overline{K}/K)$ by$$ \xymatrix{\Gal(\overline{K}/K)\ar@{>>}[r]& \Gal(F/K)\ar[r]& \GL_2(\C).}$$

\item Let $\epsilon$ be the sign representation of $S_3$: $\epsilon(s)=1$ and $\epsilon(t)={-}1$. This gives a $1$\--dimensional Artin representation  of  $\Gal(\overline{K}/K)$ which is the same Artin representation as when taking $F=\Q_5(\zeta_3)$ and  the non\--trivial $1$\--dimensional representation $\Gal(F/K)\cong C_2\to \GL_1(\C)=\C^\times$.

\item The trivial representation $\rho=\mathbf{1}$ gives the same Artin representation over~$K$ for all $F/K$. 
\end{itemize}
\end{exa}

\underline{Notation}: $\,$\\
When $K/\Q_p$ is finite and $F/K$ is a Galois extension, we write:
\begin{itemize}
\item $\pi_K$ for a fixed uniformizer of $K$;
\item $\O_K$ for the ring of integers of $K$;
\item $v_K$ for the normalized valuation on $K$;
\item $\F_K$ for the residue field of $K$ (of characteristic $p$), with $q=\#\F_K$;
\item $I_{F/K}$ for the inertia group, i.e.\ the subgroup of $\Gal(F/K)$ that acts trivially on the residue field of $F$: 
$$I_{F/K}\colon=\left\{g\in \Gal(F/K): g(x)\equiv x \bmod \pi_F \;\forall x\in \O_F\right\},$$
when $F/K$ is finite.
\item $\Frob_{F/K}$ for a Frobenius element, i.e.\ any $g\in \Gal(F/K)$ that induces the Frobenius automorphism $x\mapsto x^q$ on the residue field of $F$. 
\item $\Phi_{F/K}=\Frob_{F/K}^{{-}1}$ for the ``geometric Frobenius".
\end{itemize}

\begin{dfn} For an Artin representation $\rho$ over a local field $K$ its \emph{local polynomial} is 
$$P(\rho, T)=\det(1{-}\Phi_{F/K}T|\rho^{I_{F/K}})$$
where $\rho$ factors through $F/K$ and $\rho^{I_{F/K}}$ is the subspace of $\rho$ of $I_{F/K}$\--invariant vectors.  
\end{dfn}
Note that as $I_{F/K}$ is a normal subgroup, 
$\rho^{I_{F/K}}$ 
is a subrepresentation of $\rho$, and that $P(\rho, T)$ is essentially the characteristic polynomial of $\Phi_{F/K}$ on this subspace.
\bigskip
\bigskip

\begin{exa}
Let $K=\Q_5$ and $F=\Q_5(\zeta_3, \sqrt[3]{5})$, as in the previous example. Here
$I_{F/K}\cong C_3=\Gal(F/K(\zeta_3))$ and $\Frob_{F/K}=t$, 
an element of order $2$.\\
For the trivial representation, $P(\mathbf{1}, T)=\det(1-\triv(t)T|\triv^{I_{F/K}})=\det(1{-}T|\mathbf{1})=1-T$.\\
For the sign representation, $P(\epsilon, T)=\det(1-\epsilon(t)T|\epsilon^{I_{F/K}})=\det(1{+}T|\epsilon)=1+T$.\\
For the 2-dimensional irreducible representation of $S_3$, we find that $P(\rho, T)=1$, 
since $\rho^{I_{F/K}}=0$.
\end{exa}
\bigskip
\begin{dfn}For an {Artin representation} $\rho$ over a number field $K$, its \emph{Artin L-function} is 
$$L(\rho, s)=\prod_{\p\subset \O_K {\text{ prime}}}\frac{1}{P_\p(\rho, \nm(\p)^{-s})},$$
where $P_\p(\rho, T)$ is the local polynomial for $\rho$ restricted to $\Gal(\overline{K}_\p/K_\p)$. The Euler product is known to converge for $\Re(s)>1$ to an analytic function. 
\end{dfn}
\bigskip
\begin{exa}
If $\rho=\mathbf{1}$, then $P_\p(\rho, T)=1-T$ for all $\p$, so we obtain the Dedekind $\zeta$\--function of $K$:
$$L(\mathbf{1},s)=\prod_{\p}\frac{1}{ 1-\nm(\p)^{-s}}=\zeta_K(s).$$

Let $K=\Q$ and $\rho$ be the order $2$ character of $\Q(\zeta_3)/\Q$. Then  $P_{(3)}(\rho, T)=1$,  
$P_p(\rho, T)=1\!-\!T$ for $p\equiv 1 \pmod 3$ and  $P_p(\rho, T) =1\!+\!T$ for $p\equiv 2 \pmod 3$. So 
$$L(\rho, s)=\prod_{p\neq 3}\frac{1}{1-\left(\frac{p}{3}\right)p^{-s}}=\sum_{n=1}^\infty\left(\frac{n}{3}\right) n^{-s}=
L\left(\left(\frac{\cdot}{3}\right),s\right),$$
the $L$-function of the non-trivial Dirichlet character $\Z/3\Z\to \C^\times$.
\end{exa}
Fact: Artin $L$\--functions of $1$\--dimensional Artin representations over $\Q$ correspond to Dirichlet $L$\--functions of primitive characters.

\subsection*{Basic properties}

(i) For  $\rho_1$ and $\rho_2$ Artin representations over a local field $K$,  
we clearly have $P(\rho_1\oplus\rho_2, T)=P(\rho_1, T) P(\rho_2, T)$.

(ii) When $F/K$ is a finite extension and $\rho$ is an Artin representation over $F$,  then $P_{F}(\rho, T^{f})=P_{K}(\Ind\rho, T)$ where $f$ is the residue degree of $F/K$ and  $P_F(\rho, T)$ is the local polynomial for the Artin representation $\rho$ over $F$ and analogously for  $P_{K}(\Ind\rho, T)$.

(iii) When $K$ is a number field, $L(\rho_1\oplus\rho_2, s)=L(\rho_1, s) L(\rho_2, s)$. If $F/K$ is a finite extension and $\rho$ an Artin representation over $F$, then $L(\rho, s)=L(\Ind\rho, s)$, where the first is an Artin $L$-function over $F$ and the second over $K$.

\begin{rmk} Artin $L$\--functions of $1$\--dimensional representations are known to be analytic on $\C$ (except for a pole at $s=1$ for $\rho=\mathbf{1}$). Brauer's induction theorem taken together with (iii)  then shows that all  Artin $L$\--functions are meromorphic.  
\begin{conj}[Artin]
Let $\rho\neq \mathbf{1}$ be an irreducible Artin representation over a number field. Then its $L$\--function  is analytic.
\end{conj}
\end{rmk}

\begin{dfn}
For an Artin representation $\rho$ over a local field $K$, its \emph{conductor exponent} $n_{\rho}$ is
$$n_{\rho}=n_{\rho,\mbox{\scriptsize tame}}+n_{\rho,\mbox{\scriptsize wild}}$$
with $$n_{\rho,\mbox{\scriptsize tame}}=\dim \rho-\dim \rho^I=\dim \rho/\rho^I,$$ and  
$$n_{\rho,\mbox{\scriptsize wild}}=\sum_{k=1}^\infty \frac{1}{[I:I_k]} \dim \rho/\rho^{I_k};$$
where $\rho$ factors though $\Gal(F/K)=G$, and $I=I_{F/K}=I_0$ and
$$I_k=\{\sigma\in G: \sigma(\alpha) \equiv \alpha\bmod \pi_F^{k+1}\;\;\forall \alpha \in \O_F\}$$
are the higher ramification groups (with the lower numbering).
So, in particular, $$I_1={\rm Syl}_p I = \mbox{  wild inertia, }\;\qquad I/I_1=\mbox{ tame inertia (cyclic)}.$$
We say that $\rho$ is \emph{unramified} (resp. \emph{tame}) if $n_\rho=0$ (resp. if $n_{\rho,\mbox{\scriptsize wild}}=0$), equivalently, if $I$ (resp. $I_1$) acts trivially on $\rho$.  The \emph{conductor} of $\rho$ is the ideal $\cnd_\rho=(\pi^{n_{\rho}})$.
\end{dfn}
\begin{rmks} 
We clearly have: $n_{\rho_1 \oplus\rho_2}=n_{\rho_1}+n_{\rho_2}$ and $\cnd_{\rho_1\oplus\rho_2}=\cnd_{\rho_1}\cnd_{\rho_2}$.

One can show that when $\rho$ is irreducible and ramified,
$$n_{\rho,\mbox{\scriptsize wild}}=\dim(\rho)\cdot \max\{i: G^i \;\mbox{ acts non\--trivially on } \rho\},$$
where $G^i$ is the $i$\--th ramification group in the upper numbering.
\end{rmks}

\pagebreak

\begin{thm}[Swan's character] Let $\rho$ be an Artin representation of a local field $K$, which factors through $\Gal(F/K)$. Then 
$$n_{\rho,\mbox{\scriptsize wild}}=\langle \Tr\rho, b\rangle,$$ where 
\begin{equation*}
b(g)= \left\{\begin{array}{ll}
1-v_F(g(\pi_F)-\pi_F) & \mbox{ for } g \in I\setminus \{e\}\\
-\sum_{h\neq e}b(h)  & \mbox{ for }  g=e,
\end{array}\right.
\end{equation*}
$v_F$ is the normalized valuation of $F$ and $\langle\cdot,\cdot\rangle$ is the representation theoretic inner product for characters of $I_{F/K}$.
\end{thm}

\begin{thm}[Artin]$n_\rho\in \Z$.\end{thm}

\begin{thm}[Conductor\--discrimant formula]$\,$\\
\begin{minipage}{0.2\linewidth}
$$
\xymatrix{
F\ar^{H}@{-}[rd] \ar_{G}@{-}[dd] &  \;\\
\; & L \ar@{-}[dl] \\ 
K &\, }$$
\end{minipage}
\begin{minipage}{0.8\linewidth} 
Let $F/K$ be a Galois extension with an intermediate field $L$, and let $\rho$ be a representation of $H=\Gal(F/L)$. Then 
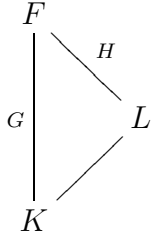
$$n_{\Ind_H^G \rho}=(\dim \rho)\cdot v_K(\Delta_{L/K})+f_{L/K} \cdot n_\rho.$$
Equivalently, 
$$\cnd_{\Ind \rho}=\Delta_{L/K}^{\dim \rho}\cdot \nm_{L/K}(\cnd_{\rho}),$$
where $f_{L/K}$ is the residue degree, $\Delta_{L/K}$ is the relative discriminant, and  $\nm_{L/K}$ the relative norm.
\end{minipage}
\end{thm}

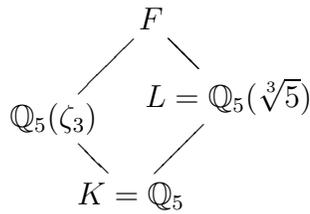
\begin{exa}$\,$\\
\begin{minipage}{0.4\linewidth}
$$
\xymatrix@dr@C=1pc{
F\ar@{-}[r] \ar@{-}[d]  &L=\Q_5(\sqrt[3]{5}) \ar@{-}[l] \ar@{-}[d] \\ 
\Q_5(\zeta_3) \ar@{-}[r] &K=\Q_5 }$$
\end{minipage}
\begin{minipage}{0.6\linewidth} 
$I_{F/K}\cong C_3=\Gal(F/K(\zeta_3))$\\
$I_1=\{1\}$\\
$\,$\\
$n_{\mathbf{1}}=0$\\
$n_\epsilon=0$\\
$n_{\rho,\mbox{\scriptsize tame}}=2-0=2, \quad n_{\rho,\mbox{\scriptsize wild}}=\sum 0=0,\quad n_\rho=2.$
\end{minipage}
In particular, by the conductor-discriminant formula,
$$1 \cdot \Delta^1_{L/K}=\cnd_{\Ind_{C_2}^{S_3}\mathbf{1}}=\cnd_\rho\cdot \cnd_{\mathbf{1}}=5^2 \;\mbox{up to units,}$$
$$\Delta_{F/K}=\cnd_{\Ind_{1}^{S_3}\mathbf{1}}=\cnd_{\rho\oplus \rho\oplus \epsilon \oplus \mathbf{1}}=5^4 \;\mbox{up to units}.$$
\end{exa}
\bigskip
\begin{dfn}
The  \emph{conductor} of an Artin representation $\rho$ over a number field $K$ is 
$$\cnd_\rho=\prod_{\p \mbox{\scriptsize{\emph{ prime in }}}K}\p^{n_\p(\rho)},$$  where $n_\p(\rho)$ is the conductor exponent of $\rho$ restricted to $\Gal(\overline{K}_\p/K_\p)$.\end{dfn}
\bigskip
\begin{thm}The Artin $L$-function of $\rho$ satisfies the functional equation:
$$\Lambda(\rho, s)=w \cdot A^{\frac{1}{2}-s} \cdot \Lambda(\rho^\ast, 1-s),$$
where
$$\Lambda(\rho, s)= L(\rho,s) \cdot \prod_{\nu \mbox{\emph{\scriptsize{ real}}}}\Gamma_\R(s)^{d_{+}(\rho)}\Gamma_\R(s+1)^{d_{-} (\rho)} \prod_{\nu \mbox{\scriptsize{\emph{ complex}}}}\Gamma_\C(s),$$ 
$\rho^\ast$ denotes the dual representation, $d_{\pm}(\rho)$ is the dimension of the $\pm$\--eigenspace of the image of complex conjugation at~$\nu$,  
$w\in \C^\ast$ with $|w|=1$ is the \emph{global root number} and 
\begin{eqnarray*}
A&=& \nm(\cnd_\rho)\cdot \sqrt{|\Delta_K|}^{\;\dim \rho},\\
\Gamma_\R(s)&=& \pi^{{-}\frac{s}{2}}\Gamma(s/2),\\
\Gamma_\C(s)&=&(2\pi)^{{-}s}\Gamma(s).
\end{eqnarray*}
\end{thm}

\newpage

\section*{Exercises to Lecture 1}
\bigskip
\setcounter{ex}{0}
\begin{ex}
For each irreducible Artin representation over $\Q$ that factors through  $\Q(\zeta_3, \sqrt[3]{3})$ determine its conductor and the first $5$ coefficients of its Artin L\--series $\sum_{n \geq 1} a_n n^{-s}$.
\end{ex}
\bigskip\bigskip\begin{ex}$\,$ \\
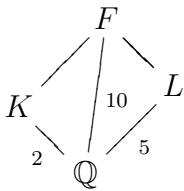
\begin{minipage}{0.2\linewidth}
$$
\xymatrix@dr@C=1pc{
F\ar@{-}[r] \ar@{-}[d] \ar@{-}^{10}[dr] & L \ar@{-}[l] \ar@{-}^{5}[d] \\ 
K \ar@{-}_ {2}[r] &\Q }$$
\end{minipage}
\begin{minipage}{0.8\linewidth} Suppose $F/\Q$ is a Galois extension with Galois group $D_{10}$, the dihedral group of order $10$. Let $K$ and $L$ denote intermediate fields with $[K:\Q]=2$ and $[L:\Q]=5$. Prove the identity $$\zeta_F(s)\zeta_\Q(s)^2=\zeta_L(s)^2\zeta_K(s).$$
\end{minipage}
\end{ex}
\bigskip\bigskip
\begin{ex}
The polynomial $f(x)=x^5+2x^4-3x^3+1$ has Galois group $D_{10}$. Its complex roots are 
$$\alpha_1=-3.01\dots, \qquad \alpha_2=-0.35 \ldots-0.53\dots i, \qquad \alpha_3=0.85\ldots-0.31\dots i,$$
$$\alpha_4=\overline{\alpha_3}, \qquad \alpha_5=\overline{\alpha_2},$$
and the Galois group contains the $5$\--cycle $(\alpha_1 \alpha_2\alpha_3\alpha_4\alpha_5)$. Prove that the Frobenius element of a prime above $2$ in the splitting field of $f(x)$ has order $5$ and determine its conjugacy class in the Galois group.\\
Hint: look at $$\alpha_1\alpha_2+\alpha_2\alpha_3+\alpha_3\alpha_4+\alpha_4\alpha_5+\alpha_5\alpha_1,$$ $$\alpha_1\alpha_3+\alpha_3\alpha_5+\alpha_5\alpha_2+\alpha_2\alpha_4+\alpha_4\alpha_1,$$and   $$\beta_1\beta_2+\beta_2\beta_3+\beta_3\beta_4+\beta_4\beta_5+\beta_5\beta_1,$$ where $\beta_1$ is a root of $f(x)$ in $\overline{\F}_2$ and $\beta_i=\beta_1^{2^{i-1}}$.
\end{ex}

\newpage 
\section{$l$-adic representations}
\begin{dfn}A continuous $l$\--adic representation over $K$ is a continuous homomorphism $\Gal(\overline{K}/K)\to \GL_d(\mathcal{F})$ for some finite extension $\mathcal{F}/\Q_l$.
\end{dfn}
\bigskip
\begin{rmk}
An $l$\--adic representation is continuous if an only if for all $n$ there exists a finite Galois extension $F_n/K$ such that   $\Gal(\overline{K}/F_n)\to \Id\; \bmod \;l^n$,
i.e.\ such that $\rho \bmod \;l^n$ factors through a finite extension $F_n/K$ (except that the image of $\rho$ may have denominators, so $\rho\bmod l^n$ may not actually be well-defined).

An Artin representation $\rho$ that factors through $F/K$ can be realized over $\overline{\Q}$, and hence also over a finite extension $\mathcal{F}/\Q_l$. It maps $\Gal(\overline{K}/F)$ to $\Id$, so we can take $F_n=F$ for all $n$ to see that it is continuous.
\end{rmk}
\bigskip
\begin{exa}[Cyclotomic character]
Let $\zeta_{l^n}$ be primitive $l^n$\--roots of unity in $\overline{K}$ 
with
$(\zeta_{l^n})^l=\zeta_{l^{n{-}1}}$. 
For $g\in \Gal(\overline{K}/K)$ define a sequence of integers $0\le a_i <l$~by 
\begin{eqnarray*}
g(\zeta_l)&=&\zeta_l^{a_1},\\
g(\zeta_{l^2})&=&\zeta_{l^2}^{a_1+l\, a_2},\\
&\dots &\\
g(\zeta_{l^n})&=&\zeta_{l^n}^{a_1+l\, a_2+\cdots+l^{n-1}a_n}.
\end{eqnarray*}
We then define the {\em $l$-adic cyclotomic character} $\chi_{cyc}$ by
$$
  \chi_{cyc}(g)=a_1+l\, a_2+\cdots+l^{n-1}a_n+\cdots\in\Z_l.
$$
Note that the value $\chi_{cyc}(g)\bmod l^n$ simply says what $g$ does to the $l^n$\--th roots of $1$. It is easy to check that the $l$\--adic cyclotomic character is multiplicative, and hence gives a 1-dimensional representation,
$$\chi_{cyc}: \Gal(\overline{K}/K)\to \Z_l^\ast\subset\GL_1(\Q_l).$$
Taking $F_n=K(\zeta_{l^n})$, we have that $\Gal(\overline{K}/F_n)\to \Id \bmod\, l^n$, so $\chi_{cyc}$ is continuous.
\end{exa}
\bigskip

\begin{exa}[Tate module]
Let $E/K$ be an elliptic curve 
and let $P_n, Q_n$ be a basis for $E[l^n]$ with $lP_n=P_{n-1}$ and $lQ_n=Q_{n-1}$. 
For $g\in \Gal(\overline{K}/K)$, we define $0\le a_i, b_i, c_i, d_i<l$ by
$g(P_1)= a_1 P_1 + b_1 Q_1$ and $g(Q_1)=c_1 P_1+ d_1 Q_1$, and generally
$$
\begin{array}{cccccccccc}
g P_n &=& (a_1 + \ldots + a_nl^{n-1}) P_n &+ &(b_1 + \ldots + b_nl^{n-1}) Q_n, \cr
g Q_n &=& (c_1+\ldots+c_nl^{n-1}) P_n &+ & (d_1+\ldots +d_nl^{n-1}) Q_n.
\end{array}
$$
Then 
$$
\rho(g)=\abc{a_1+\ldots+l^{n-1}a_n+\ldots\quad}{c_1+\cdots+l^{n-1}c_n+\cdots}{b_1+\ldots+l^{n-1}b_n+\ldots\quad}{d_1+\ldots+l^{n-1}d_n+\ldots} 
\in \GL_2(\Z_l)\subseteq \GL_2(\Q_l)
$$ 
is the representation on the $l$-adic Tate module of $E$. Here we can take $F_n=K(E[l^n])$ to see that $\rho$ is continuous, and the value of $\rho(g) \bmod l^n$ says what $g$ does to $E[l^n]$.
\end{exa}
\newpage
\begin{rmks}$\,$

\begin{enumerate}[(i)]
\item Let $K/\Q_p$ be a finite extension and $\mathcal{F}/\Q_l$ with $l\neq p$. Then $\Frob_{\overline{K}/K}(\zeta_{l^n})=\zeta_{l^n}^q$ and $I_{\overline{K}/K}$ acts trivially on $\zeta_{l^n}$ for all $n$. In other words $$\chi_{cyc}(I_{\overline{K}/K})=1 \qquad {\text{and}} \qquad \chi_{cyc}(\Frob_{\overline{K}/K})=\#\F_K=q.$$
\item If $\rho$ is the representation on the Tate module of an elliptic curve, then $\det \rho=\chi_{cyc}$, i.e.\ the determinant of the matrix $\rho(g)$ simply depends on the action  of  $g$ on $\zeta_{l^n}$ for all $n$. This is because the Weil pairing $E[l^n]\times E[l^n]\to\mu_{l^n}$ is alternating and Galois equivariant (and is, effectively, the determinant map).
\item For elliptic curves, when $K/\Q_p$ is a finite extension and $l\neq p$, one usually uses the dual representation $$\rho_E=\rho^\ast,$$ i.e.\ $\rho_E(g)=(\rho(g)^{-1})^t$; this is the same as the representation on $\coh^1_{\mbox{\scriptsize{\emph{\'et}}}}(E, \Q_\ell)$.
\end{enumerate}
\end{rmks}\bigskip
\begin{dfn}
Let $K/\Q_p$ be a finite extension  and $\rho:\Gal(\overline{K}/K)\to \GL_d(\mathcal{F})$ a continuous $l$\--adic representation with $l\neq p$. 
The \emph{local polynomial} of $\rho$ is $$P(\rho, T)=\det(1{-}\Phi_{\overline{K}/K} T|\rho^{I_{\overline{K}/K}}),$$
and its {\em conductor exponent} is $n_{\rho}=n_{\rho,\mbox{\scriptsize tame}}+n_{\rho,\mbox{\scriptsize wild}}$ with 
\begin{eqnarray*}
n_{\rho,\mbox{\scriptsize tame}}&=&\dim \rho/\rho^{I_{\overline{K}/K} },\\
n_{\rho,\mbox{\scriptsize wild}}&=&\sum_{k\geq 1} \frac{1}{[I_{F/K}:I_{F/K, k}]} \dim \rho/\rho^{I_{F/K, k}},
\end{eqnarray*}
where $F/K$ is a finite extension, large enough so that the action of the wild inertia group factors through $F/K$. 
(This exists: one can take $F=F_1$ so that the image of $\Gal(\overline{K}/F)$ lies in 
$\abcd{1{+}l\mathcal{O_F}}{l\mathcal{O_F}}{l\mathcal{O_F}}{1{+}l\mathcal{O_F}}$.
The image of wild inertia is trivial since it is a (pro) $p$\--group and this matrix group a (pro) $l$\--group). The {\em conductor} of the representation $\rho$ is $\cnd_\rho=(\pi_K)^{n_{\rho}}$.
\end{dfn}
\begin{dfn}For an elliptic curve $E$ over a number field $K$, its \emph{L-function} is defined by the Euler product
$$L(E, s)=\prod_{\p \mbox{ \scriptsize{\emph{prime}}}}\frac{1}{P_\p(\rho_E, \nm(\p)^{-s})},$$
where $P_\p(\rho_E, T)$ is the local polynomial for $\rho_E:\Gal(\overline{K}_\p/K_\p)\to \GL_2(\Q_l)$ for any  $l$ not divisible by $\p$. It converges for $\Re(s)\gg 1$. Its conductor is $\cnd_E=\prod_\p \p^{n_{\p,\rho_E}}$ and it conjecturally satisfies the functional equation:
$$\Lambda(\rho, s)=w \cdot A^{1-s}\cdot  \Lambda(\rho, 2-s),$$
where
\begin{eqnarray*}
\Lambda(\rho, s)&=&L(E,s) \Gamma_\C(s)^{[K:\Q]},\\
\Gamma_\C(s)&=&(2\pi)^{{-}s}\Gamma(s),\\
A&=& \nm_{K/\Q}(\cnd_E)\cdot |\Delta_K|^2,\\
w&=&\pm 1 \mbox{ is the \emph{global root number} of }E/K .
\end{eqnarray*}
\end{dfn}

\begin{rmk}This construction of $L$\--functions applies more generally to  ``compatible systems of $l$-adic representations " (e.g.\ from Artin representations, abelian varieties,  $\coh^i_{\mbox{\scriptsize{\emph{\'et}}}}$ of varieties, motives).
The exact formula for the $\Gamma$\--factors and the root number $w$ is explicitly given using the associated Hodge structure and the theory of local $\epsilon$\--factors of $l$-adic representations.
\end{rmk}

\section{The $l$-adic representation of an elliptic curve}

Recall that to an elliptic curve $E/K$ we have associated an $l$-adic representation $\rho_E$. We now relate its properties to the arithmetic of $E$ in the case when $K$ is a local field.

\begin{thm}
Let $K/\Q_p$ be a finite extension, $E/K$ an elliptic curve, and let \hbox{$\rho_E:\Gal(\overline{K}/K)\to \GL_2(\Q_l)$} the $l$-adic representation on $\coh^1_{\mbox{\scriptsize{\emph{\'et}}}}(E, \Q_\ell)$ with $l\neq p$. 
Then
\begin{enumerate}[1.]
\item $E$ has good reduction if and only if $\rho_E$ is unramified (the N\'eron\--Ogg\--Shafarevich criterion),
\item $\det\rho_E=\chi_{cyc}^{-1}$ (by the Weil pairing),
\item $P(\rho_E, {1}/{q})={\#\tilde{E}(\F_K)}/{q}$, where $q=\#\F_K$.
\end{enumerate}
\end{thm}
\bigskip
\begin{rmks}
\begin{enumerate}[(i)]
\item By $(1)$ $E$ has potentially good reduction if and only if $I_{\overline{K}/K}$ acts through a finite quotient.
\item By $(1), (2)$ and $(3)$, if $E/K$ has good reduction then 
$$P(\rho_E, T)=1-aT+qT^2\quad\mbox{with}\quad a=1+q-\#\tilde{E}(\F_K).$$
\item By $(1)$ and $(3)$:
\begin{eqnarray*}
E \mbox{ has additive reduction } &\Rightarrow & P(\rho_E, T)=1,\\
E \mbox{ has split multiplicative reduction } &\Rightarrow & P(\rho_E, T)=1-T,\\
E \mbox{ has non\--split multiplicative reduction } &\Rightarrow & P(\rho_E, T)=1+T.
\end{eqnarray*}
\end{enumerate}
\end{rmks}
\newpage
\begin{exa}$\,$

\begin{minipage}{0.6\linewidth}
Let $E/\Q_5$ be the elliptic curve $y^2=x^3+5^2$. It has additive reduction ($\Delta_E= \mbox{unit} \cdot 5^4$). 
It has good reduction over  $L=\Q_5(\sqrt[3]{5})$: over $L$ we have the model $E': y^2=x^3+1$ with $\tilde{E}'(\F_5)=\{\O, (0,\pm 1),(2, \pm2),({-}1,0) \}$. 
The action of $\Gal(\overline{\Q}_5/\Q_5)$ factors through $\Gal(L^{nr}/\Q_5)$ which is generated by  $g\in \Gal(L^{nr}/\Q^{nr}_5)$ of order $3$ (inertia group) and 
$\Phi_L=\Phi_{L^{nr}/L}$, the geometric Frobenius over $L$. They satisfy the relation  $$\Phi_L \cdot g \cdot \Phi_L^{-1}=g^{-1}.$$
\end{minipage}
\begin{minipage}{0.4\linewidth}
$$
\xymatrix@dr@C=1pc{
\overline{\Q}_5 \ar@{-}[r] & L^{nr} \ar^{C_3=\langle g \rangle}@{-}[r] \ar@{-}[d] & \Q^{nr}_5\ar@{-}[dd]  \\ 
\; & L(\zeta_3)\ar@{-}[d] \ar^{S_3}@{.}[dr] &\;\\
\; & L\ar^3 @{-}[r]\ar^{\langle \Phi_L \rangle}@/^2pc/@{-}[uu]  &\Q_5}$$
\end{minipage}
 In particular $\rho_E(g)$ has eigenvalues $\zeta_3, \zeta_3^{-1}$ (it is nontrivial as there is bad reduction, so it has order $3$, and it has determinant $1$ by the Weil pairing) and $\rho_E(\Phi_L)$ has eigenvalues $\pm \sqrt{{-}5}$ (as $\#\tilde{E}'(\F_5)=6$ so $P(\rho_E, T)=1+5T^2$ over $L$); and $$\rho_E(\Phi_L) \cdot \rho_E( g) \cdot \rho_E(\Phi_L)^{-1}=\rho_E(g)^{-1}.$$ 

A little algebra shows that, with a suitable choice of $\overline{\Q}_l$\--basis, the representation $\rho_E: \Gal(L^{nr}/K)\to \GL_2(\overline{\Q}_l)$ is given by 
$$\rho_E(g)=\abcd{\zeta_3}{0}{0}{\zeta_3^{-1}} \quad\quad\quad \mbox{ and } \quad\quad\quad\rho_E(\Phi_L)=\abcd{0}{\sqrt{{-}5}}{\sqrt{{-}5}}{0}.$$
\end{exa}

\newpage 

\section*{Exercises to Lecture 2}
\bigskip
\setcounter{ex}{0}
\begin{ex}
Let $E/\Q_7$ be the elliptic curve $y^2+y=x^3-x^2$. Let $K/\Q_7$ be some horrible extension with residue degree $11$ and ramification degree $76$. 
Find $\#\tilde{E}(\F_K)$, where $\F_K$ is the residue field of $K$.
\end{ex}
\bigskip
\bigskip
\begin{ex}
Let $K/\Q_p$ be a finite extension and $E/K$ an elliptic curve. 
Prove that if $p\neq 2,3$ then the exponent of the conductor of $E$ is at most $2$.

Hint: pick $\ell\neq p$ that forces the image of the wild inertia group in $\GL_2(\Z/\ell^n\Z)$ to be trivial for all $n$.
\end{ex}
\bigskip
\bigskip
\begin{ex}
Let $E/\Q_7$ be the elliptic curve $y^2=x^3+7^2$. Identify its $\ell$\--adic representation 
$\rho_E:\Gal(\overline{\Q}_7/\Q_7)\to \GL_2(\overline{\Q}_\ell)$ on 
$\coh^1_{\mbox{\scriptsize{\emph{\'et}}}}(E, \Q_\ell)\otimes_{\Q_\ell}\overline{\Q}_\ell$ for $l\neq 7$. 

Hint: compute $\#\tilde{E}(\F_K)$ over two different cubic ramified extensions $K/\Q_7$.
\end{ex}

\newpage 
\section{Classification of $l$-adic representations}
Throughout this section $K$ is a finite extension of $\Q_p$ and $l\neq p$ a prime.

\begin{exa}
$\,$\\

\begin{minipage}{0.6\linewidth}
In the last lecture we had $E/\Q_5$ with $l$\--adic representation $\rho_E$ given by 
$\rho_E(g)=\abcd{\zeta_3}{0}{0}{\zeta_3^{-1}}$ and $\rho_E(\Phi_{L})=\abcd{0}{\sqrt{{-}5}}{\sqrt{{-}5}}{0}$, 
where $\langle g \rangle =I\cong C_3$ with $g( \sqrt[3]{5})= \zeta_3 \sqrt[3]{5}$, and $\Phi_{L}=\Phi_{L^{nr}/L}$ is a topological generator of  $\Gal(L^{nr}/L)$.

The representation $\rho_E$  can clearly be written as 
$$\rho_E=\rho\otimes \psi,$$ 
where $\rho$ is the $2$\--dimensional representation of $\Gal(\Q_5(\zeta_3, \sqrt[3]{5})/\Q_5)\cong S_3$: 
$$\rho(g)=\abcd{\zeta_3}{0}{0}{\zeta_3^{-1}}, \quad \rho_E(\Phi_{L})=\abcd{0}{1}{1}{0},$$
\end{minipage}
\begin{minipage}{0.4\linewidth}
\xymatrix@dr@C=1pc{
\overline{\Q}_5 \ar@{-}[rr] & \, & \Q^{nr}_5(\sqrt[3]{5}) \ar@{-}[r]^(0.65){C_3=\langle g \rangle} \ar@{-}[d] & \Q^{nr}_5\ar@{-}[d]  \\ 
\; &\; &  L(\zeta_3) \ar@{-}[d] \ar@{.}[dr]^(0.35){S_3}\ar@{-}[r] & \Q_5(\zeta_3)\ar@{-}[d] \\
\; & \; &L{=}\Q_5(\sqrt[3]{5})\ar^{\langle \Phi_{L} \rangle}@/^2pc/@{-}[uu]  \ar @{-}[r]  &\Q_5}
\end{minipage}
and $\psi$ satisfies $\psi(I)=1$ and $\psi(\Phi_{\overline{\Q}_5/\Q_5})={\sqrt{{-}5}}$. 

This decomposition works as $\rho_E(\Phi_L)^2=\abcd{-5}{0}{0}{-5}$ is a scalar matrix. 
Indeed, $\Phi_L^2$ is central in $\Gal( \Q^{nr}_5(\sqrt[3]{5})/\Q^{nr}_5)$ and $\rho_E$ is irreducible, so by Schur's lemma $\rho_E(\Phi_L)^2$ must be scalar. 
This trick works in general for irreducible representations, and lets  one prove the following classification:
\end{exa}
\bigskip
\begin{thm} Every continuous $l$\--adic representation $\tau: \Gal(\overline{K}/K)\to \GL_d(\mathcal{F})$ for which 
\begin{itemize}
\item the image of $I$ is finite,
\item any ($\Leftrightarrow$ every) choice of $\Phi_{\overline{K}/K}$ acts semisimply,
\end{itemize}
is of the form $$\tau \cong\bigoplus_{i}\rho_i\otimes \chi_i\quad\quad\quad\quad\mbox{ (for $\mathcal{F}$ sufficiently large),}$$
with
\begin{itemize}
\item $\rho_i$ Artin representations (with values in $\mathcal{F}$),
\item $\chi_i$ $1$\--dimensional unramified: that is $\chi_i(I)=1$ and $\chi_i(\Phi_{\overline{K}/K})\in \mathcal{F}$.
\end{itemize}
\end{thm}
It is not hard to show that for irreducible $l$\--adic representations with finite image of inertia, $\Phi_{\overline{K}/K}$ must act semisimply.

\newpage

\begin{exa}
$\,$\\

\begin{minipage}{0.7\linewidth}
Let us now consider an $l$\--adic representation for which the image of inertia $I=I_{\overline{K}/K}$ is infinite. Let $E$ be an elliptic curve over $K$, with $K/\Q_p$ finite, such that $E$ has split multiplicative reduction,  and let 
$\rho_E:\Gal(\overline{K}/K)\to \GL_2(\Z_l)$ be the associated $l$-adic representation. Then
$P(\rho_E, T)=1-T$, so there is a $1$\--dimensional $I$\--invariant space with trivial $\Phi_{\overline{K}/K}$ action:
$$\rho_E(h)=\abcd{1}{?}{0}{?}\quad  \forall h\in \Gal(\overline{K}/K).$$ 
Since $\det\rho_E=\chi_{cyc}^{-1}$, it follows that $\det\rho_E(h)=1$ for all $h\in I$, and so $\rho_E(h)=\abcd{1}{\ast}{0}{1}$ for $h\in I$ with $\ast \in \Z_l$ and not always zero (as otherwise there would exist a $2$\--dimensional inertia invariant subspace, which  is not the case). 
\end{minipage}
\begin{minipage}{0.3\linewidth}
\xymatrix@dr@C=1pc{
\overline{\Q}_p \ar@/^2pc/@{-}[rrr]^{I}\ar@{-}[rr]& \,  & L^{nr} \ar@{-}[r]^(0.3){\langle g\rangle}\ar@{-}[d]^{\Phi_L}& K^{nr}\ar@{-}[d] \\ 
\; & \;& L\ar@{-}[r]& K}
\end{minipage}
Hence, the image of inertia mod $l^n$ is of the form
$$
\rho_E(I)\subseteq\left\{\abcd{1}{\ast}{0}{1}\in \GL_2(\Z/l^n\Z)\right\}\cong C_{l^n}.
$$
As for every $n$ it is contained in a cyclic group of order $l^n$, the wild inertia (a $p$-group, $p\neq l$) maps to $\Id$, and the tame inertia 
$\Gal({K}^{\mbox{\scriptsize{tame}}}/K^{nr})\cong \prod_{p'\neq p}\Z_{p'}$ is such that only $\Z_l$ can act non\--trivially (and must do so). 
So, if $\Z_l=\langle g\rangle=\Gal(L/K^{nr})$, where $L=\cup_{n=1}^\infty K^{nr}({\sqrt[l^n]{\pi_K}})$ and $\Phi_L=\Phi_{L^{nr}/L}$, then $\rho_E$ factors through $\Gal(L^{nr}/K)$ and, with respect to a suitable $\Q_l$\--basis, 
$$
\rho_E(g)=\abcd{1}{1}{0}{1}\quad\quad \quad \mbox{and }\quad\quad\quad \rho_E(\Phi_L)=\abcd{1}{0}{0}{q}$$
where $q=\#\F_L$. More explicitly,
$$\rho_E(h)=\abcd{1}{t_l(h)}{0}{\chi^{-1}_{cyc}(h)}\quad\mbox{for }h\in I,$$where $t_l:I\to\Z_l$ is the $l$\--adic tame character given by $h(\sqrt[l^n]{\pi_K})= \zeta_{l^n}^{t_l(h)}\sqrt[l^n]{\pi_K}$ for all $n$.
\end{exa}

\begin{rmk}
Once we have $\rho_E(g)=\abcd{1}{t_l}{0}{1}$ and $\rho_E(\Phi_l)=\abcd{1}{?}{0}{?}$, the rest is, in fact, forced by the commutation relations of $\Phi_L$ and $I$. In particular, $\det\rho_E$ has to be $\chi_{cyc}^{-1}$. 
\end{rmk}

\begin{dfn} The \emph{special representation} $sp(n)$ over $K$ is the $n$\--dimensional $l$-adic representation given by:
$$sp(n)(h)=
\left(\begin{matrix} 
1 &  t  & t^2/2  & \cdots & t^{n{-}1}/(n{-}1)!\\
 0& 1  & t & \cdots & t^{n{-}2}/(n{-}2)!\\
 & &  \ddots& \ddots & \vdots\\
 & & & & t\\
 0 & & & 0 & 1
\end{matrix}\right) \quad\mbox{ for }h\in I,$$
where $t=t_l(h)$ is the $l$\--adic tame character, and

$$sp(n)(\Phi_{\overline{K}/K})= \left(\begin{smallmatrix} 
1 &   &  & 0  \\
 & {q} &  &   \\
 & &  \ddots&  \\
 0 & & & q^{n{-}1}\end{smallmatrix}\right),\quad\mbox{ where }q=\#\F_K.$$
(Different choices of $\Phi_{\overline{K}/K}$ give rise to isomorphic representations.)

Note that, in particular, $sp(1)$ is the trivial representation, and $sp(2)$ agrees with the one associated to an elliptic curve with split multiplicative reduction.
\end{dfn}

\begin{thm} Every continuous $l$\--adic representation $\tau: \Gal(\overline{K}/K)\to \GL_d(\mathcal{F})$ for which  any ($\Leftrightarrow$ every) choice of $\Phi_{\overline{K}/K}$ acts semisimply on $\tau^{I'}$, for every  $I'\subseteq I$ of finite index,  is of the form 
$$\tau=\bigoplus_{i}\rho_i \otimes sp(n_i)\quad\quad\mbox{ (for $\mathcal{F}$ sufficiently large),}$$
where $n_i$ are integers and $\rho_i$ are continuous $l$\--adic representations with $\rho_i(I_{\overline{K}/K})$ finite and semisimple action of $\Phi_{\overline{K}/K}$.
\end{thm}
\begin{exa}
Let $E/K$ be an elliptic curve with split multiplicative reduction. Then $\rho_E$ has infinite inertia image, so, by the theorem, $\rho_E=\rho\otimes sp(2)$ where $\rho$ is $1$\--dimensional with local polynomial $1-T$. Therefore, $\rho={\mathbf{1}}$ and $\rho_E(g)=sp(2)$.
\end{exa}
\bigskip

\section{Independence of $l$}

Throughout this section $K$ is a finite extension of $\Q_p$ and $l\neq p$ a prime.

\begin{exa}
The $l$\--adic cyclotomic character is defined such that $\chi_{cyc}(I)=1$ and $\chi_{cyc}(\Phi_{\overline{K}/K})=1/q$, where $q=\#\F_K$. Morally, this does not depend on $l$. However, $\Gal(\overline{K}/K)$ is a topological group, and limits do depend on $l$. For instance, the $\chi_{cyc}$ factors through $\cup_{n}K(\zeta_{l^n})/K$, which depends on $l$. 
\end{exa}
\medskip
\begin{exa}
In the last lecture we had $E/\Q_5$ with $\rho_E(g)=\abcd{\zeta_3}{0}{0}{\zeta_3^{-1}}$ and  $\rho_E(\Phi_L)=\abcd{0}{\sqrt{{-}5}}{\sqrt{{-}5}}{0}$. 
This representation does not depend on $l$ provided we only look at the inertia group and {\em integer} powers of $\Phi_L$. Non-integer powers of $\Phi_L$ again give problems because of different convergence properties in $\Q_l$ for different $l$.
\end{exa}
\medskip
\begin{exa}
Let $E/K$ be an elliptic curve with split multiplicative reduction. Then $\rho_E(g)=sp(1)$ which again, morally, does not depend on $l$. (However, the $l$\--adic tame character does).
\end{exa}
\medskip
\begin{dfn}
The \emph{Weil group} $W_{\overline{K}/K}$ is the subgroup of $\Gal(\overline{K}/K)$ consisting of elements whose image modulo $I_{\overline{K}/K}$ is an integer power of $\Phi_{\overline{K}/K}$, i.e.\
$$\xymatrix{
0\ar[r]&I_{\overline{K}/K}\ar[r]\ar@{=}[d]& \Gal(\overline{K}/K)\ar[r]&\widehat{\Z}\cong\Gal(\F_{\overline{K}}/\F_K)\ar[r]&0\\
0\ar[r]&I_{\overline{K}/K}\ar[r]&\;W_{\overline{K}/K}\;\ar[r] \ar@{}_{\rotatebox{90}{$\subseteq$}}[u]&\Z\cong\langle \Phi_{\overline{K}/K}\rangle\ar@{}^{\rotatebox{90}{$\subseteq$}}[u] \ar[r] &0.
}$$
The topology of $W_{\overline{K}/K}$ is the same (profinite) one on $I_{\overline{K}/K}$ and discrete on $W_{\overline{K}/K}/I_{\overline{K}/K}$.
\end{dfn}
\begin{rmk}
A continuous $l$\--adic representation over $K$ automatically restricts to a continuous representation $W_{\overline{K}/K}\to \GL_d(\mathcal{F})$.
\end{rmk}
\medskip
\begin{thm}
Let $E/K$ be an elliptic curve (or an abelian variety). Then the decomposition of $\rho_E$ as $\bigoplus_i \rho_i\otimes sp(n_i)$ is independent of $l$, i.e.\ $\rho_i\otimes_{\Q_l}\C$ are independent of $l$ as representations of $W_{\overline{K}/K}$, and the $n_i$ are independent of $l$.
\end{thm}

\begin{cor}
Let $E/K$ be an elliptic curve (or an abelian variety) with potentially good reduction. Then $\rho_E\otimes_{\Q_l}\C$ is independent of $l$ as a representation of $W_{\overline{K}/K}$.
\end{cor}
\medskip


\begin{cor}
Let $E/K$ be an elliptic curve (or an abelian variety). The local polynomial $P(\rho_E,T)$ and the characteristic polynomials of $\rho_E(h)$ for every $h\in I_{\overline{K}/K}$ 
are independent of $l$ and lie in $\Q[x]$ (as they lie in $\Q_l[x]$ for all $l\neq p$).
\end{cor}

\medskip
\begin{exa}Let $K/\Q_p$ be a finite extension with $p\geq 5$, and let $E/K$ be an elliptic curve with potentially good reduction. Then $\rho_E$ is tame (exercise $2$ of lecture $2$), 
so the action of $I_{\overline{K}/K}$ factors through a finite cyclic group $C_n=\langle g \rangle$ (as the tame inertia group is cyclic). 
By the above corollary, we must have $C_n=C_1, C_2,C_3, C_4$ or $C_6$ since the characteristic polynomial of $\rho_E(g)$ is quadratic and lies in $\Q[x]$.
 (Thus, in particular, $v_K(\Delta_E)\not\equiv 1,11\bmod 12$ and $E$ acquires good reduction over a quartic or a sextic totally ramified extension of $K$).
\end{exa}

\newpage
\section*{Exercises to Lecture 3}
\bigskip
\setcounter{ex}{0}
\begin{ex}
Show that if $\rho_E$ is irreducible then $E$ acquires good supersingular reduction over a finite extension. Recall that $E/K$ has good supersingular reduction if  $\#\tilde{E}(\F_K)\equiv 1 \bmod p$.
\end{ex}
\bigskip
\bigskip
\begin{ex}Show that an elliptic curve with potentially multiplicative reduction acquires split multiplicative reduction over a quadratic extension.
\end{ex}
\bigskip
\bigskip
\begin{ex}Decompose $sp(2)\otimes sp(2)$ as a direct sum of indecomposables. (If you like the representation theory of $\GL_2(\C)$ try also $sp(n)\otimes sp(m)$). 
\end{ex}

\end{document}